\documentclass[letterpaper, 10 pt, conference]{ieeeconf}  

\IEEEoverridecommandlockouts                              

\usepackage{graphicx} 
\usepackage{epsfig} 
\usepackage{mathptmx} 
\usepackage{times} 
\usepackage{amsmath} 
\usepackage{amssymb}  
\usepackage{color}
\usepackage{subfigure}
\usepackage{psfrag}

\newcommand{\bma}{\begin{bmatrix}}
\newcommand{\ebma}{\end{bmatrix}}
\newcommand{\bd}{\mathbf}
\newtheorem{theorem}{Theorem}
\newtheorem{corollary}{Corollary}

\title{\LARGE \bf
Stability of Interconnected DC Converters
}

\author{Gustavo Cezar$^{1}$, Ram Rajagopal$^{2}$ and Baosen Zhang$^{3}$
\thanks{$^{1}$Gustavo Cezar is a graduate student with the Department of Civil and Environmental Engineering at Stanford University
        {\tt\small gcezar@stanford.edu}.}%
\thanks{$^{2}$Ram Rajagopal is an Assistant Professor in Civil and Environmental Engineering and, by courtesy, Electrical Engineering and Center Fellow at the Precourt Institute for Energy at Stanford University {\tt\small ramr@stanford.edu}.}%
\thanks{$^{3}$Baosen Zhang is an Assistant Professor in Electrical Engineering at University of Washington
        {\tt\small baosen.zhang@gmail.com}.}%
\thanks{This work is partially supported by the Tomkat Center for Sustainable Energy, E-ON, China EPRI and a Google Award.} 
}

\begin{document}

\maketitle
\thispagestyle{empty}
\pagestyle{empty}

\newcommand{\todo}[1]{\textcolor{red}{GC: #1}}
\newcommand{\todoB}[1]{\textcolor{red}{BZ: #1}}

\begin{abstract}

This paper addresses stability issues of DC networks with constant power loads (CPL). Common DC networks, such as automotive electrical systems and DC microgrids, typically have a step-up/down converter connected in one side to the main bus and, on the other, to the load. When load is constant power it can generate destabilizing effects if not proper controlled.  This paper shows that converters driving CPLs can make the system unstable, even if they are individually stable, depending on network  parameters. We mitigate this problem by means of passive components externally connected to the converter/CPL subsystem. The analysis is verified through simulations. We are able to show that certain converter circuit configurations achieve the so called plug-and-play property, which stabilizes the interconnected system for all network parameters. This property is desirable since it is does not require the knowledge of detailed system topology and parameters, which can be time varying and difficult to obtain. This method also contrasts to existing practices of load augmentation, which can lead to severe efficiency losses. 

\end{abstract}

\section{Introduction}
The electricity grid is undergoing a significant change in both generation and demand sides. Increasing penetration of renewables has led to increased deployment of wind and solar generation, in particular distributed photovoltaics.  The profile of loads has also faced a significant change with increased adoption of electronic devices such as computer power supplies, LED lights, cellphones electric vehicles and home storage units. These different supply and demand technologies are connected to the electric grid through power electronics converters forming an interconnected network. 

Examples of power electronic components are DC/AC (or AC/DC) inverter and DC/DC converter. Inverters are used to transform alternating current (AC) systems to direct current (DC) systems, and vice-versa. Converters are usually used to step-up or step-down DC voltages to meet bus and loads operating points. Traditionally, inverters and converters have been studied as individual components, focusing on their efficiency and other performance criteria. 

A main challenge in design of such components is ensuring stable operations. These devices present a destabilization characteristic due to their intrinsic nonlinearity when tightly regulated in the presence of loads that are constant power \cite{emadi2}. Constant power loads (CPL) display a negative impedance characteristic \cite{grig}\cite{belk}\cite{emadi3}. Operations instability arises as the converter tries to keep the output power constant. Thus, when the input current decreases, the input voltage increases (and vice-versa). The aggregate effect is of a negative impedance that acts as a positive feedback to the system.  Prior research has paid significant attention to this issue and studied the stability of a single source, inverter and load system based on linear, non-linear and passivity based control techniques \cite{emadi1,emadi2,emadi3,alexis1,alexis2,amr}.  \emph{Yet, not much is known in terms of power electronics control stability in networks  of converters connected to constant power loads}. 

In this paper, we focus on a network of connected DC/DC converters. Such a network could arise in practice as a DC Microgrid, a standalone power system on large ships and planes, or even a home with multiple types of distributed energy resources. We choose to focus on DC networks\footnote{This is different then the DC approximation used in for AC power systems. Here by DC system we mean that all lines are resistive and voltages are not sinusoidal.} because most new devices such as rooftop solar and storage are inherently DC devices.  

We show that even though single converters may be stable, the network is not necessarily stable. In fact, we show that if each converter only has access to its own states, then no feedback controller can ensure that the system is stable for all network parameters. Also, we show that as the size of the network increases, the system becomes unstable for all most all line resistance parameters. Intuitively, the disturbances in the network propagate from converter to converter and is not controllable even if each converter is stable by themselves. 

We design different external circuit topologies where passive components, such as capacitors and resistors, are added in different parts of the system. We show that with input capacitive filters, we can achieve a plug-and-play property that ensures that a network is stable for all line resistance parameters. Our method does not require modifications to the internal structure of converters, allowing off-the-shelf components to be used.  Front-end filter contrasts with existing designs in \cite{jusoh,passdamp}, which consider passive damping strategies by inserting passive components either to the filter elements in the converter or in parallel with the load at the output of the converter. 

The paper is organized as follows.  In Section II, the concepts of CPLs and negative impedance are presented. Also, a derivation of the dynamics of a buck converter operating in continuous conduction mode feeding CPL is shown. In Section III a simple nonlinear control is develop to make a single buck converter feeding a CPL stable. Also we show that instability may arise when cascading converters on the same bus and show through simulations that the more converters are added to the bus, more likely it is for the system to become unstable. Section IV studies different scenarios where passive components can be added to stabilize the system. Finally, concluding remarks are presented in Section V.



\section{Converter and Load Models} \label{sec:model}

The increased number of sources and electronic loads that need to regulate their voltage in order to meet bus voltage and load operating conditions creates stability issues in DC systems. The  issues arise as voltage regulation occurs by means of converters. Converters are devices that regulate voltage but keep power constant. When a converter is connected to a load, it forms a system that can be represented as a CPL. Ideally, a CPL can be represented by:

\begin{equation}\label{CPL}
    I_o(t)=\frac{P}{V_o(t)}, \; \forall\;  V_o^{\min} \leq V_o(t) \leq V_o^{\max}.
\end{equation}

\noindent where $P$ is the power required by the load, $V_o(t)$ is the load operating voltage, and $I_o(t)$ is the current draw by the load. Constant power loads behavior is qualitatively different than more standard load models such as constant resistance loads. As voltage increases, constant power loads draw less current, whereas constant resistance loads draw more current. The latter acts as a damping element for a circuit, whereas the former acts as a positive feedback and tends to amplify disturbances in the system. In practice, this leads to large voltage and current swings, and the load shuts off if voltage is outside of the design range $[V_o^{\min},V_o^{\max}]$. Next, we describe the mathematical model used to analyze buck converters driving constant power loads. 

\subsection{Buck Converter Modeling}
We concentrate on the ideal buck converter in this paper due to its popularity. Similar analysis can be carried out for other types of DC/DC converters.
\begin{figure}[thpb]
      \centering
      \includegraphics[scale=0.7]{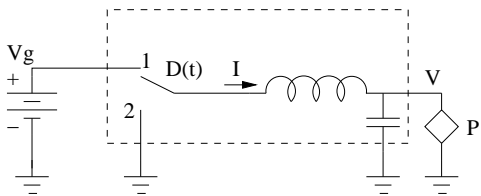}
      \caption{Buck converter feeding a constant power load. $D(t)$ is the duty cycle of the switch, which is the main control knob of the system.}
      \label{fig:buck_two}
   \end{figure}
   
An ideal Buck converter (see Fig. \ref{fig:buck_two}) is a switched circuit, where the switching duty cycle is denoted $D(t)$ by power electronic convention \cite{Kas}. The notation is slightly misleading since $D(t)$ is not a function of time, rather it is a fraction of nominal period $T$. Rigorously, let $t|T$ denote the modulo operation, that is, $t|T=t-kT$ where $k$ is the unique integer where $kT \leq t < (k+1)T$. The duty cycle $D(t)$ is defined only for $t=nT$ for integer $n$ and takes on values between $0$ and $1$. Then the switch is at position $1$ (connected to the main line) if $\frac{t|T}T \leq D(kT)$; otherwise the switch is at position $2$ (grounded). Therefore to be rigorous, the system needs to be thought as a discrete time system modulo $T$. However, most converters have switching frequency of tens to hundreds of kHz, so $T$ is much shorter than any other time constants/periods in the system. Hence we can think of the circuit as a switched continuous time system where $D(t)$ is understood to mean the (possibly changing) duty cycle. 

The state of the Buck converter can be described by the current through the inductor and the voltage across the capacitor, denoted as $I(t)$ and $V(t)$, respectively. The input to the control system is changes in the duty cycle.  The evolution of the state depend on the position of the switch:
\begin{itemize}
\item Switch at position 1 (connected to main line):
\begin{equation} \label{eqn:1}
L \frac{d I}{d t} = V_g - V(t) \quad C \frac{d V}{d t} = I(t)-\frac{P}{V(t)}.
\end{equation}
\item Switch at position 2 (grounded):
\begin{equation} \label{eqn:2}
L \frac{d I}{d t} = - V(t) \quad C \frac{d V}{d t} = I(t)-\frac{P}{V(t)}.
\end{equation}
\end{itemize}
The \emph{average signal model} enables a continuous time analysis of switched electronic circuits that can summarize eqs.~\eqref{eqn:1} and \eqref{eqn:2}.   
In this model, for a given signal $x(t)$, its average is defined by
\begin{equation}
<x(t)>= \frac{1}{T} \int_{t}^{t+T} x(\tau) d \tau,
\end{equation}
where $T$ is the some time period that is often taken to be the switching period.  By applying the average model to the converters equations \eqref{eqn:1} and \eqref{eqn:2}, the states of the Buck converter can be expressed as:
\begin{equation} \label{eqn:Indaverage}
L \frac{d <I>}{d t} = D(t) (<V_g> - <V>),
\end{equation}

\begin{equation} \label{eqn:Capaverage}
C \frac{d <V>}{d t} = <I>-\frac{P}{V}.
\end{equation}

Eqs~\eqref{eqn:Indaverage} and \eqref{eqn:Capaverage} form a set of differential equations neither linear nor time invariant. In order to asses stability and evaluate controllers a LTI model of the circuit needs to be obtained. A common technique used is the \textit{small signal analysis}. By assuming that currents and voltages of a converter consist of a dc component and a small ac variation, the average model can be linearized by Taylor expansion. For an average signal $<x(t)>=X+ x(t)$, where $X$ is the dc component and lower case $x(t)$ is the small ac variation. Therefore, \eqref{eqn:Indaverage} and \eqref{eqn:Capaverage} can be rewritten as:

\begin{equation} \label{eq:IndSm}
L \frac{d (I+ i(t))}{d t} = [D+d(t)] [V_g+  v_g(t)] - [V+ v(t)], \\
\end{equation}

\begin{equation} \label{eq:CapSm}
C \frac{d (V+ v(t))}{d t} = (I+ i(t))-\frac{P}{(V+ v(t))}.
\end{equation}

Expanding \eqref{eq:IndSm} and \eqref{eq:CapSm} and approximating the nonlinear term $P/(V+ v(t))$ to $P/V - v(t)P/V^2$, making the derivatives of the DC component on the left side equal to zero, neglecting second order terms, and canceling the DC components on the right side, the small signal LTI model is written as:

\begin{equation} \label{eq:IndSmA}
L \frac{d ( i(t))}{d t} = [D v_g(t)+V_g d(t)]- v(t) \\
\end{equation}

\begin{equation} \label{eq:CapSmA}
C \frac{d (\hat v(t))}{d t} =  i(t)-\frac{ v(t) P}{V^2}
\end{equation}

Letting $x=[ v(t) \quad \hat i(t)]$ and $u= d(t)$. Therefore the duty cycle can be thought of as the control knob in the system. Writing \eqref{eq:IndSmA} and \eqref{eq:CapSmA} in matrix form as:

\begin{equation} \label{eq:ConMat}
\dot{\textbf{x}}=\textbf{Ax}+\textbf{B}\textit{u}
\end{equation}

\noindent
where

\begin{equation} \label{eqn:Mat}
\textbf{A}= \bma \frac{P}{C \bar{V}^2} & \frac{1}{C} \\ -\frac{1}{L} & 0 \ebma \;\; \textbf{B}=\bma 0 \\ \frac{V_g}{L} \ebma,
\end{equation}
and $\bar{V}$ is the average steady state value of the output voltage.

\section{Feedback Control} \label{sec:stability}

\subsection{Feedback Control of Single Converter}
To analyze the dynamics of a Buck converter feeding a CPL one has to look at the eigenvalues of the matrix \textbf{A} in \eqref{eq:ConMat}. These eigenvalues are given by:

\begin{equation}
\lambda = \frac{1}{2}\{ P/C\bar{V}^2 \pm \sqrt{(P/C\bar{V}^2)^2-4/LC} \}.
\end{equation}

Therefore the system with a CPL is open loop unstable due to a pair of eigenvalues with positive real part. However, it can be easily shown that this system is controllable, thus there exists a closed loop controller.

\begin{figure}[thpb]
      \centering
      \includegraphics[scale=0.5]{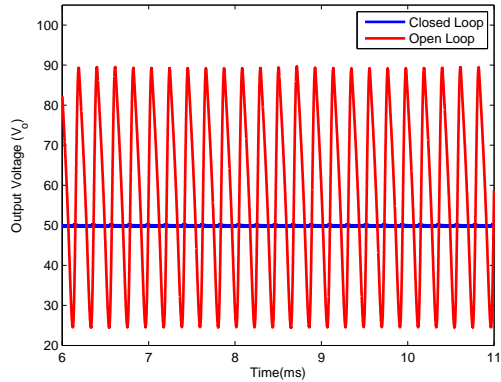}
      \caption{Open loop and closed loop output voltage response of a buck converter driving a CPL}
      \label{fig:OLCL}
   \end{figure}

To show the open loop instability behavior of a Buck converter and the closed loop response, a simulation was performed with the standard software \textit{PSIM}. The parameters of the circuit were: $V_i=110VDC$, $L=20\mu H$, $C=29\mu F$, $P_o=1kW$, and $f_{sw}=200kHz$. Fig. \ref{fig:OLCL} shows open loop and closed loop output voltage. They are overlayed in the same graph for comparison. For the open loop case, the system response is oscillatory with output voltage swings around 60V. This result, in practice, will never happen because loads operate only in a narrow band of voltages determined by their specifications. In the closed loop case, using a proportional controller, output voltage swings significantly decreased, with a peak to peak amplitude of around 0.2V.

\subsection{Stability of Cascaded Converters}
\begin{figure}[thpb]
      \centering
      \includegraphics[scale=0.7]{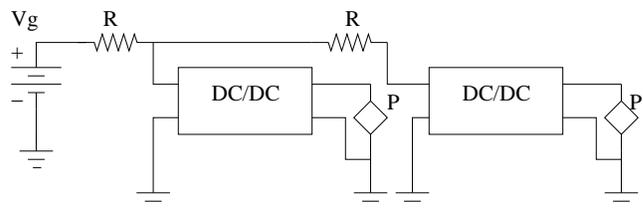}
      \caption{Cascaded converters driving constant loads. The lines are modeled as resistors with value $R$.}
      \label{fig:CascConver}
   \end{figure}

The main focus of design has been on a single converter because the performance of many applications such as automotive drives are dominated by the main converter. With increased penetration of renewables and electronic based loads, converters are likely to be connected in a \emph{network}. For example, Fig. \ref{fig:CascConver} shows two cascaded converters. We say that a converter is individually stable if it is stable by itself as a standalone circuit connecting a power source to a constant power load.

Theorem \ref{thm:two_converter} shows that no matter how the feedback control is designed for each individual converter, there are networks parameters where the whole system is not stable.
\begin{theorem}\label{thm:two_converter}
Suppose two Buck converters are connected in series as in Fig. \ref{fig:CascConver}, each driving a constant load. Assume that all connecting lines have the same resistance $R$. Suppose each converters is individually stabilized by linear state feedback control. Then for all linear feedback controllers, there are values of $R$ such that the network is not stable.
\end{theorem}
 The proof of Theorem \ref{thm:two_converter} is given in the Appendix and Fig. \ref{fig:two_converter} illustrates the unstable behavior. The equal resistance assumption simplifies the derivation but is not required.  Most commercial converters are embedded with a feedback control mechanism to ensure that it is stable when driving constant power loads. Theorem \ref{thm:two_converter} shows that no matter how individual controllers are designed, the networked system can be unstable.

\begin{figure}[thpb]
      \centering
      \includegraphics[scale=0.6]{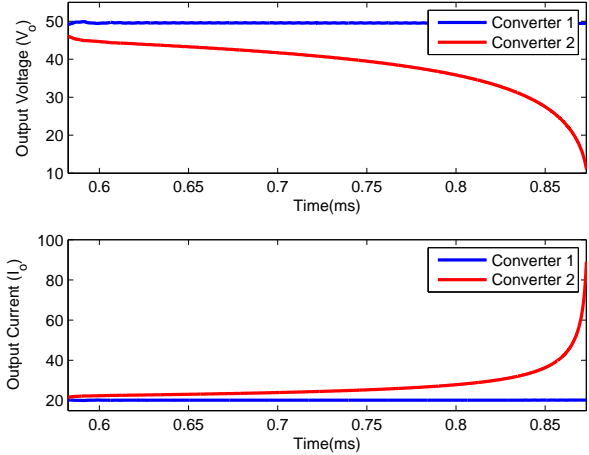}
      \caption{Instability of cascaded converters}
      \label{fig:two_converter}
\end{figure}

\begin{figure}[thpb]
      \centering
      \includegraphics[scale=0.5]{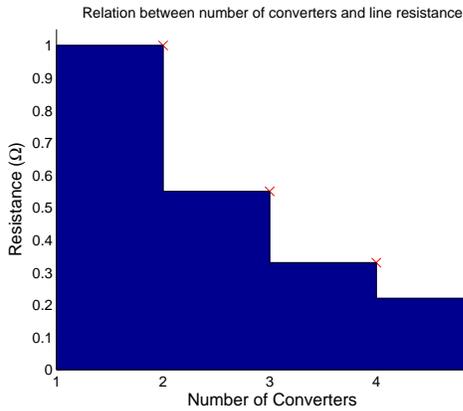}
      \caption{Relation between number of converters and line resistances. The shaded region is where the system is stable.}
      \label{fig:many_converter}
\end{figure}

Figure \ref{fig:many_converter} shows  the situation worsens as the number of converters in the network increases. It plots the largest line resistance such that the cascaded system is stable where each converter is individually stable. As the number of converters increases, the maximum line resistance decreases to $0$. Since controllers are often hardcoded in commercial converters, this means that a large network of converters is likely to be unstable. Corollary \ref{cor:many_converter} formalizes this phenomenon.

\begin{corollary} \label{cor:many_converter}
Consider a line network with $n$ converters driving constant power loads. Suppose all the lines have resistance $R$. Fix feedback controllers of each individual converter. For any values of $R >0$, there exists $N_0$ such that the network is unstable if $n > N_0$.
\end{corollary}

\subsection{Distributed Control}
We have demonstrated that feedback stability of individual converters does not guarantee stability of a network converters.  One possible solution to this problem is to globally design a controller based on the output of all converters. However, since converters switch at frequencies of tens to hundreds of kHz and have no communication capabilities, sharing information is impractical. Therefore, the feedback controller is sparsity constrained, where the duty cycle of a converter in the network only depend on its own measurements. This problem is a challenging due to its non-convexity, and have recently received a lot of attention (see, e.g. \cite{Fazelnia14, DJ13} and references within).  However, this approach is still hard to implement since it requires the knowledge of network topology and parameters of all components. In practice, devices are connected in ad-hoc basis and converters from different manufactures have different parameters. In the next section, we design converters with plug-and-play capabilities \emph{by modifying passive electronics}.

\section{Cascaded converter stability by means of passive components}
The last section showed that network converters are unstable and designing distributed controllers require topology and parameter information that is often not available. In this section, we explore ways to ensure that converters have plug-and-play capabilities. That is, the networked system is guaranteed to be stable for all parameters.

\begin{figure}[ht]
\centering
\psfrag{Vin}{$V_{in}$}
\subfigure[Resistance at output]{\label{fig:out_R}
\includegraphics[scale=0.6]{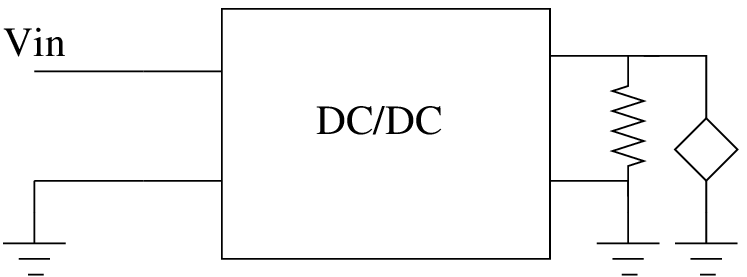}}
\subfigure[Input RC filter at ground]{\label{fig:in_RC}
\includegraphics[scale=0.6]{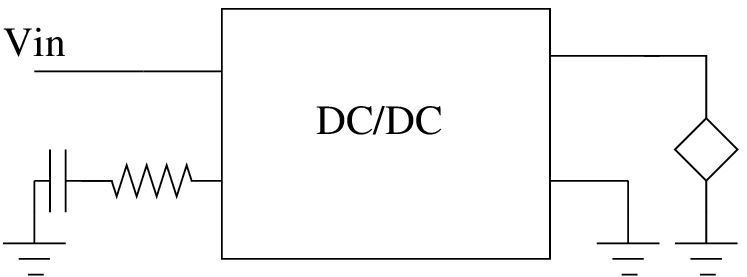}}
\subfigure[Input shunt capacitor]{\label{fig:input_C}
\includegraphics[scale=0.6]{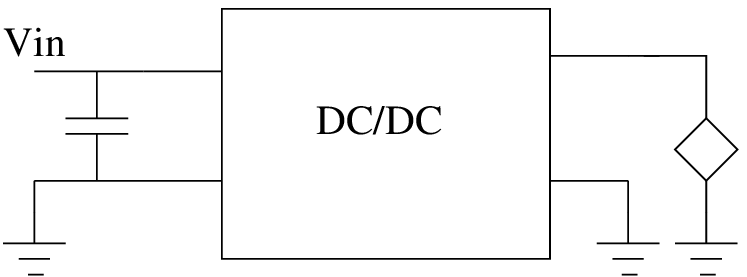}}
\caption{Three possible designs to change the behavior of the converter. The top figure shows the most popular design with a shunt resistance at the output. The middle figure adds a RC filter between ground and the converter. The bottom figure adds a shunt capacitor at the input. Out of the three, we show the bottom design is the best one.}
\label{fig:designs}
\end{figure}

Due to space constraints and for simplicity, we focus on a network with two cascaded converters as in Fig. \ref{fig:CascConver}. We assume that the internal structure of the converters cannot be changed and any modification must be achieved in the circuitry around the converters. Fig.~\ref{fig:designs} shows three possible designs to change the behavior of the converters.
The performance capabilities of these designs are stated in Theorem \ref{thm:design}.
\begin{theorem}\label{thm:design}
\textit{Consider a network of two cascade converters where each one is modified according to one of the circuits in Fig. \ref{fig:designs}. If both converters are modified as Fig. \ref{fig:input_C}, then the network is stable for all line resistances for some capacitor value. If both converters are modified as Fig. \ref{fig:out_R}, to guarantee stability, at least half of the total power would be dissipated in the shunt resistor. If the both converters are modified as Fig. \ref{fig:in_RC}, the system cannot be guaranteed to be stable.}
\end{theorem}
The proof is stated in the Appendix.

Out of the three designs in Figure \ref{fig:designs}, the top design is the most popular. It adds a shunt resistor to counteract the effects of the constant power load \cite{emadi2}. Intuitively, for a small enough shunt resistance, the aggregated load at the output of the converter acts like a resistor, which eliminates instabilities in the system. However, as stated in Theorem \ref{thm:design}, to achieve the desired effect, the power lost in the resistor is at least as much as the power consumed by the load. Therefore the efficiency of converter is at most $50\%$, which is not acceptable in most applications.
The middle design is has the desired properties that the RC filter does not effect other components in the system, but it cannot make the network stable.

\begin{figure}[thpb]
      \centering
      \includegraphics[scale=0.6]{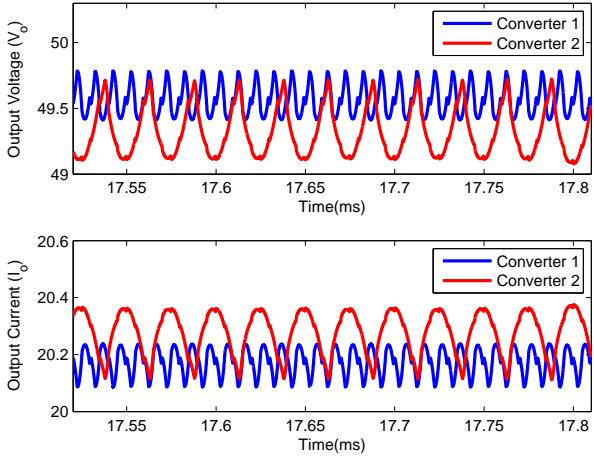}
      \caption{Output voltages and currents for schematic of Fig. \ref{fig:input_C}}
      \label{fig:input_C_plot}
\end{figure}

The bottom design in Fig. \ref{fig:designs} can stabilize the system for all line resistances. In essence, the capacitor acts as a storage device that damps out the current propagation between two cascaded converters. Therefore disturbances are localized, and feedback controllers on each individual converter are able to stabilize the system. Fig. \ref{fig:input_C_plot} plots the states of the system.

Figure \ref{fig:input_C_plot} shows that for a system with individually stable converters driving a CPL but unstable when connected in a network (refer to Fig. \ref{fig:CascConver}), can be made stable by adding a shunt capacitor to the converter's input. This shunt capacitor, as mentioned before, will act as a storage device in order to minimize the current drawn from the power source. Without the shunt capacitor, whenever the current on the constant power load starts increasing, more current will be drawn from the power source. With this higher current coming from the source, the voltage drop across the line resistances will increase. Thus, the voltage at the input of the converter will decrease and ultimately the load current will increase, generating a destabilization effect. The shunt capacitor will prevent more current being draw from the source by providing additional current to the system. Thus, this destabilization effect is canceled. In addition, this solution comes at a lower cost if compared to other solutions such as placing a resistor in parallel with the load. The losses due to the equivalent series resistance (ESR) of the capacitor are considerably smaller than the losses in the load's parallel resistor. First because the ESR is a fraction of the parallel resistor needed, and second because the shunt capacitor does not need to provide all the load current, but just a small variation.

\textbf{Note:} The proof of Theorem \ref{thm:design} relies on the linear approximation of the system, which leads to input filters with large $C$. Very large capacitors could introduce unwanted ripples into the actual nonlinear system. Better understanding performance and stability trade-offs in the nonlinear system is an important area for future work.

\section{Conclusions}
This paper analyzed stability issues of dc systems. Most practical dc systems can be viewed as a variation of this problem where a main bus is connected to loads through a step-down/up converter, i.e. automotive electrical system and DC microgrids. These loads are seeing by the converter as a constant power and thus introduce destabilizing effects on the system such as voltage oscillations and collapse.

It was shown that even if a converter driving a CPL can be made individually stable in closed loop, when connected in a network, it can drive the system to an unstable condition. Theorem \ref{thm:two_converter} shows that independently of the feedback controller used in each individual converter, there will be always network parameters that drive the system to an unstable condition. To solve this problem, instead of designing a robust controller or adding components in parallel with the load, which most of the time is not feasible to implement in commercially available products, this paper takes a different approach. It addresses the stability issue externally to the converter by means of passive components. It was shown that adding a shunt capacitor in parallel with each converter, a system that was previously unstable can be made stable. Future work will focus on hardware implementation and testing and also study the destabilization effects of CPL on AC systems.
\appendices
\section{Proof of Theorem \ref{thm:two_converter}}
For simplicity of notation, we assume that both converter have the same values for $C$,$L$ and $P$. Fix feedback controllers for each of the converters. We show that as the line resistance $R$ increases, the overall system is not stable. 

The state of the system is $[I_1 \; V_1 \; I_2 \; V_2]^T$ where $I_1$ and $I_2$ are the currents through inductors in converter 1 and 2, respectively; $V_1$ and $V_2$ are the voltage across the capacitors in converter 1 and 2, respectively. As in Section \ref{sec:model}, we use $<>$ to denote average quantities and lower cased letters to denote small signal quantities. Since averaging is a linear operation, when considering the switched operation of converter $1$, the signals from converter 2 is taken to be its averaged value (similarly for converter 2). When the switch of converter 1 is connected to the line, the state evolution equations are
\begin{equation*}
L \frac{d I_1}{dt}=V_g - (I_1+<I_2>) R-V_1 \quad C \frac{d V_1}{dt} = I_1-\frac{P}{V_1},
\end{equation*} 
if the switch is grounded, 
\begin{equation*}
L \frac{d I_1}{dt}=-V_1 \quad C \frac{d V_1}{dt} = I_1-\frac{P}{V_1}. 
\end{equation*}
For converter 2, if the switch is connected to the line, the state evolution equations are 
\begin{equation*}
L \frac{d I_2}{dt}=V_g - (<I_1>+I_2) R-I_2 R -V_2 \quad C \frac{d V_2}{dt} = I_2-\frac{P}{V_2},
\end{equation*} 
if the switch is grounded, 
\begin{equation*}
L \frac{d I_2}{dt}=-V_2 \quad C \frac{d V_2}{dt} = I_2-\frac{P}{V_2}. 
\end{equation*}
Averaging and making the smaller signal assumptions, we have the following linearized system equations
\begin{equation}
\dot{\mathbf{x}}=\bd{A} \bd{x} + \bd{B} \bd u,
\end{equation}
where $\bd x=\bma i_1 & v_1 & i_2 & v_2 \ebma^T$, $\bd u = \bma d_1 & d_2 \ebma^T$ 
\begin{equation*}
\bd A=  \bma -\frac{DR}{L} & -\frac{1}{L} & -\frac{DR}{L} & 0 \\ \frac{1}{C} & \frac{P}{C\bar{V}_1} & 0 & 0 \\ -\frac{DR}{L} & 0 & -\frac{2DR}{L} & -\frac{1}{L} \\
 0 & 0 & \frac{1}{C} & \frac{P}{C\bar{V}_2} \ebma=\bma \bd A_1 & \bd A_{12} \\ \bd A_{21} & \bd A_2 \ebma,
 \end{equation*}
 \begin{equation*} 
 \bd B= \bma \frac{\bar{V}_1}{L} & \vdots & 0 \\ 0 &  \vdots & 0 \\ \dots & & \dots\\ 0 & \vdots & \frac{\bar{V}_2}{L} \\ 0 & \vdots & 0 \ebma=\bma \bd B_1 & 0 \\ 0 & \bd B_2\ebma,
\end{equation*}
and $\bar{V}_1$ and $\bar{V}_2$ are the average steady state values for $V_1$ and $V_2$, respectively. 

A linear state feedback controller is given in the form of $\bd u=\bd F \bd x$, where $\bd F$ is restricted to be in the form of 
\begin{equation*}
\bd F= \bma \bd F_1 & 0 \quad 0 \\ 0 \quad 0 & \bd F_2 \ebma
\end{equation*}
since each converter controller only has access to its own states. We say that each converter is individually stable if $\bd{A}_1+\bd{B}_1 \bd {F}_1$ is stable. Let $\tilde{\bd A}=\bd{A}+\bd{B} \bd{F}$. Because of the structure of the feedback, $\tilde{\bd A}$ has the following sign structure
\begin{equation} \label{eqn:global_feedback}
\tilde{\bd A}=\bma - & - & -\frac{DR}{L} & 0 \\ + & + & 0 & 0 \\ -\frac{DR}{L} & 0 & - & - \\
 0 & 0 & + & + \ebma=\bma \tilde{\bd A}_1 & \tilde{\bd A}_{12} \\ \tilde{\bd A}_{21} & \tilde{\bd A}_2 \ebma,
\end{equation}

We show that as $R$ increases, at least one eigenvalue of $\tilde{\bd A}$ is  positive. This is equivalent to proving that the determinant of $\tilde{\bd A}$ is negative. Following a standard result for the determinant of block matrices \cite{Shores07}, we have
\begin{equation}
\det(\tilde{\bd A})=\det(\tilde{\bd A}_1)\det(\tilde{\bd A}_2-\tilde{\bd A}_{12} \tilde{\bd A}_1^{-1} \tilde{\bd A}_{12}).
\end{equation} 
Since both eigenvalues of $\tilde{\bd A}_1$ have negative real parts (converter 1 is individually stable), $\det(\tilde{\bd A})$ is positive. By straightforward calculation, the matrix in the second term has the form
\begin{equation}
\tilde{\bd A}_2-\tilde{\bd A}_{12} \tilde{\bd A}_1^{-1} \tilde{\bd A}_{12}= \bma - & - \\ + & + \ebma -R^2 \bma c & 0 \\ 0 & 0\ebma, 
\end{equation}
where $c$ is a positive constant. Therefore for large enough $R$, the determinant of $\tilde{\bd A}_2-\tilde{\bd A}_{12} \tilde{\bd A}_1^{-1} \tilde{\bd A}_{12}$ is negative. Implying that at least one eigenvalue of $\tilde{\bd A}$ has positive real part.

\section{Proof of Corollary \ref{cor:many_converter}}
This proof follows from the proof of Theorem \ref{thm:two_converter} by repeatedly using the block determinant formula for a $n$-converter network. The calculations are standard but somewhat tedious and we skip them here due to space constraints. The final result of the calculation is that the sign of the determinate of the global system matrix is determined by the sign of 
\begin{equation*}
\det(\bma - & - \\ + & + \ebma -nR^2 \bma c & 0 \\ 0 & 0\ebma),
\end{equation*} 
where $c$ is a positive constant. As $n$ grows, the above determinate is negative, implying that not all eigenvalues have positive real parts (there are always $2n$ eigenvalues). 
\section{Proof of Theorem \ref{thm:design}}
We prove the three cases in this theorem one by one. For the circuit in Fig. \ref{fig:out_R}, following the notation in Appendix I, state matrices become 
\begin{equation*}
\bd A_1= \bma -\frac{DR}{L} & -\frac{1}{L} \\ \frac{1}{C} & \frac{P}{C\bar{V}_1}-\frac{1}{C R_s} \ebma \; \; \bd A_2= \bma -\frac{2DR}{L} & -\frac{1}{L} \\ \frac{1}{C} & \frac{P}{C\bar{V}_2}-\frac{1}{C R_s} \ebma,
\end{equation*}
where $R_s$ is the value of the shunt resistor. Due to the sparsity constraints in the feedback controller, the only the top row of $\bd A_1$ and $\bd A_2$ would be affected. Therefore $R_s$ need to be small enough such that $\frac{P}{C\bar{V}_1}-\frac{1}{C R_s}$ and $\frac{P}{C\bar{V}_2}-\frac{1}{C R_s}$ are negative, otherwise the result from Appendix I holds. The average loss in converter 1 is given by 
\begin{equation*}
\frac{\bar{V}_1^2}{R_s} \geq \bar{V}_1^2 \frac{P}{\bar{V}_1^2} = P. 
\end{equation*}
Similarly, the loss in converter is also greater than $P$. Hence at least 50\% of power is wasted in the shunt elements. 

For the circuit in Fig. \ref{fig:in_RC}, the RC filter at the input does not change the sign of the second row of the $\bd A_1$ and $\bd A_2$. That is, both matrices still have signs of $\bma - & - \\ + & + \ebma$. Repeating the calculations in Appendix I shows that this setup cannot stabilize the network. 

For the circuit in Fig. \ref{fig:input_C}, we expand the state space of to include the voltage across the shunt capacitors at the inputs to the converters. Denote them as $v_{c1}$ and $v_{c2}$ , respectively. Let the new state be $\bd x= \bma v_{c1} & i_1 & v_1 & v_{c2} & i_2 & v_2 \ebma^T$. Let $C_s$ be the value of the shunt capacitor. Following same steps as in Appendix I, the new state evolution matrices are 
\begin{equation*}
\bd A=\bma -\frac{2}{C_s R} & - \frac{1}{C_s} & 0 & \frac{1}{C_s R} & 0 & 0 \\
\frac{D}{L} & 0 & -\frac{1}{L} & 0 & 0 & 0 \\ 0 & \frac{1}{C} & \frac{P}{C \bar{V}_1^2} & 0 & 0 & 0\\
   \frac{1}{C_s R} & 0 & 0 &  -\frac{1}{C_s R} & - \frac{1}{C_s} & 0 \\ 0 & 0 & 0 & \frac{D}{L} & 0 & -\frac{1}{L} \\ 0 & 0 & 0 & 0 & \frac{1}{C} & \frac{P}{C \bar{V}_2^2} \ebma 
\end{equation*}
and 
\begin{equation*}
\bd B=\bma 0 & 0 \\ \bar{V}_{c1} & 0 \\ 0 & 0 \\ 0 & 0 \\ 0 & \bar{V}_{c2} \\ 0 & 0 \ebma.   
\end{equation*}
As $C_s$ increases, the two system becomes essentially decoupled, and the feedback in each converter is able to stabilize the system since each converter is individually stable. 

It is important to note that the analysis here depend on linearization and averaging. Therefore the adverse ripple effect of increasing $C_s$ is not evident. A direction of future work is to include the effects of a large input capacitor into the calculations.

\addtolength{\textheight}{-12cm}   

\bibliographystyle{IEEEtran}
\bibliography{Bibliograp}

\end{document}